\documentclass[12 pt]{article}
\usepackage{graphicx}
\usepackage{harpoon}
\usepackage{color}
 \usepackage{stackrel} 
\usepackage{latexsym}
\usepackage{amssymb}
\usepackage{amsmath}
\usepackage{mdframed}
\usepackage{hyperref, url}
\usepackage[show]{ed}
\usepackage{subfigure}
\usepackage[ntheorem]{empheq} 
\usepackage{amsthm, amssymb, amsfonts, latexsym}
\mathchardef\mhyphen="2D

\usepackage{pst-node}

\makeatletter
\DeclareRobustCommand{\cev}[1]{%
  \mathpalette\do@cev{#1}%
}
\newcommand{\do@cev}[2]{%
  \fix@cev{#1}{+}%
  \reflectbox{$\m@th#1\vec{\reflectbox{$\fix@cev{#1}{-}\m@th#1#2\fix@cev{#1}{+}$}}$}%
  \fix@cev{#1}{-}%
}
\newcommand{\fix@cev}[2]{%
  \ifx#1\displaystyle
    \mkern#20mu
  \else
    \ifx#1\textstyle
      \mkern#20mu
    \else
      \ifx#1\scriptstyle
        \mkern#26mu
      \else
        \mkern#26mu
      \fi
    \fi
  \fi
}

\makeatother

\newcommand{\fig}[1]{\epsfbox}
\newcommand{\bp}{\begin{minipage}{3.1cm}}
\newcommand{\ep}{\end{minipage}}
\usepackage{pbox}

\newtheorem{Definition}{Definition}[]
\newtheorem{Theorem}{Theorem}[]
\newtheorem{Lemma}{Lemma}[]
\newtheorem{Proposition}{Proposition}[]

\newtheorem{Remark}{Remark}[]
\newtheorem{Example}{Example}[]

\newtheorem*{Theorem*}{Theorem}

\makeatletter 
\renewcommand{\thefigure}{\@arabic\c@figure}
\makeatother
\usepackage{xr}

\begin{document}

\title{Embedding Information onto a Dynamical System}

\author{G. Manjunath \\
{\small \it Department of Mathematics \& Applied Mathematics, University of Pretoria, Pretoria 0028} \\ {\small Email: manjunath.gandhi@up.ac.za } \\}
\date{}
\date{}
\maketitle

"This is an author-created, un-copyedited version of an article accepted for publication/published in Nonlinearity. IOP Publishing Ltd is not responsible for any errors or omissions in this version of the manuscript or any version derived from it. 

\begin{abstract}
The celebrated Takens' embedding theorem concerns embedding an attractor of a dynamical system in a Euclidean space of appropriate dimension through a generic delay-observation map.   The embedding also establishes a topological conjugacy.  In this paper, we show how an arbitrary sequence can be mapped into another space as an attractive solution of a nonautonomous dynamical system. Such mapping also entails a topological conjugacy and an embedding between the sequence and the attractive solution spaces.   This result is not a generalization of Takens embedding theorem but
helps us understand what exactly is required by discrete-time state space models widely used in applications to embed an external stimulus onto its solution space. Our results settle another basic problem concerning the perturbation of an autonomous dynamical system. We describe what exactly happens to the dynamics when exogenous noise perturbs continuously a local irreducible attracting set (such as a stable fixed point) of a discrete-time autonomous dynamical system.    
 \end{abstract}
 
 \vspace{-0.5cm}

   \vspace{-0.5cm}
 \section{Introduction} \label{Sec_Intro}

 Scientists have been hard-pressed to explain how a cell or an organism could exhibit the necessary responsiveness and precision while responding to external stimuli. Long before scientists could explain such cryptic phenomena, engineers had developed control systems that keep our automated systems on course.
Researchers have now been endeavoring to build systems that mimic some of the human brain's functionalities. A system of that kind to respond precisely should capture the sequential or temporal information presented to it without distortion.  
In this work, a system's ability to contract its states, a notion that can be made mathematically accurate, is found to be the fundamental principle behind capturing sequential information without distortion in many natural and artificial systems.
 
Discrete-time state space models of the form \begin{eqnarray}
	x_{n+1} &=& g(u_n,x_n) \label{eqn1} \\
	y_{n+1} & =& f(x_{n+1}), \label{eqn2}
\end{eqnarray}
where $n\in \mathbb{Z}$, $u_n$ belongs to an input space $U$ and state $x_n$ belonging to a state space $X$, $g: U \times X \to X$, and a measurement $f: X \to Y$  have 
been used with great success in control theory applications in industrial and scientific applications mainly owing to the rich theory in the particular case of linear systems (e.g. \cite{kalman1962canonical,weiss2010lectures}). The dynamics obtained from \eqref{eqn1} is not limited to modelling of artificial systems, but also relevant to computation in natural systems, for instance a neutrophil, a single white blood cell chasing a bacterium \cite[movie recorded by David Rogers]{phillips2012physical} can be modeled as a single cell updating or computing its current position $x_n$ in the next time-step to $x_{n+1}$ after it receives a stimulus $u_n$ from the bacterium. In another context, the dynamics in \eqref{eqn1} can also be viewed as that of a noise-perturbed autonomous system.

Newly, inspired by computation in the human brain, the model in \eqref{eqn1}--\eqref{eqn2} with some select properties on the map $g$ have been popularised for information processing tasks in the field of reservoir computing since the heavily cited works of  \cite{jaeger2004harnessing}, and \cite{maass2002real}. Such processing of data involves mapping of the sequential information onto usually a higher-dimensional space, and it is empirically found that their representations in the higher dimensional space (e.g. \cite{jaeger2004harnessing,jaeger2001echo}) have certain exaggerated features. One of the other great advantages of mapping temporal data is that it can be mapped as an attractor of a nonautonomous system (e.g., \cite{kloeden2004uniform,manjunath2014dynamics}), and such attractivity gives numerical stability while processing 
(e.g., \cite{jaeger2004harnessing,grigoryeva2018echo,grigoryeva2019differentiable, 
manjunath2021universal,Manju_2020}).  The reservoir computing framework 
is also suited for multitasking that accomplishes information processing tasks such as predicting, filtering, classification of sequential input or data streams. 
 The multitasking is attributed to the possibility of designing different task-specific read-outs $f$ in \eqref{eqn2} while the hardware that implements \eqref{eqn1} remains task-independent. When such multi-tasking applications are required, the map $g$ in \eqref{eqn1} needs to satisfy an asymptotic state contraction property \cite{jaeger2004harnessing},  and a rigorous treatment as to why multitasking is indeed possible can be found in \cite[Theorem II.4]{cuchiero2021discrete}. This paper aims to analyse such mapping of sequential data through the dynamics defined by  $g$. 
  The findings in this paper cast new and sharp light on how such sequential information can be mapped or embedded onto another space. Here and throughout we say a topological space $X$ is embedded in another space $Y$ if there exists a function $f: X \to Y$  such that $f: X \to f(X)$ is a homeomorphism.

For every bi-infinite sequence $\bar{u}= \{u_n\}_{n\in \mathbb{Z}}$, the map $g$, which we call a driven system henceforth generates a non-autonomous dynamical system defined by a family of self-maps $\{g(u_n,\cdot)\}_{n\in \mathbb{Z}}$ on $X$. We say a driven system $g$ has the unique solution property (USP) if for each input $\bar{u}$ there exists exactly one solution which is a bi-infinite sequence $\{x_n\}$ that satisfies \eqref{eqn1}. In other words, $g$ has the USP if there exists a well-defined solution-map $\Psi$ between the infinite product space $\overline{U}$
 and $\overline{X}$ (here $\overline{Y}$ denotes $\prod_{i=-\infty}^\infty Y_i$ when $Y_i\equiv Y$)
  so that $\Psi(\bar{u})$ denotes the unique solution obtained from the input $\bar{u}$. We note that the phrase ``unique" in the USP does not imply that the solution-map  $\Psi$ is injective. For instance, for the driven system, $g(u,x)=ux/2$,  where $U=[0,1]$ and $X=[0,1]$, $g$ has the USP since  $\Psi (\bar{u})$ is identically equal to the sequence with only $0$'s. The USP property 
  in the neural network literature is equivalent to the echo state property (e.g., \cite{Manju_ESP}). We prefer to use the unique solution property terminology since the results here are not confined to the literature on neural networks.  We note that when $g$ has the USP, and $X$ has more than a single point, it would follow that $g(u,\cdot): X\to X$ is not surjective for any $u \in U$ (see Section~\ref{Sec_Prelim})  and the non-autonomous dynamics, in general, does not approach a single-set, but a sequence of sets -- that can be described through a certain topological contraction  (see Remark~\ref{Rem_Contraction}).  We also note that the USP property depends on the input space in addition to the mapping $g$.  We use the USP as the central notion of interest and connect it to other notions of the uniform attraction property, causal mapping, and causal embedding that we introduce.

Each nonautonomous dynamical system $\{g(u_n,\cdot)\}_{n\in \mathbb{Z}}$ induced by $\bar{u}$ can be rewritten (see Section~\ref{Sec_Prelim}) in a process formulation  -- the process is a two-parameter semigroup formalism (e.g., \cite{kloeden2011nonautonomous,kloeden2013discrete})
\begin{align} \label{eq_Process}
  \phi_{\bar{u}}(n,m,x) := \begin{cases}
        x \quad \quad \quad \quad \quad \quad \quad \text\quad \quad \quad \quad \quad \quad \:\text{ if $n=m$,}
        \\
        g_{u_{n-1}}\circ \cdots \circ g_{u_{m+1}} \circ g_{u_{m}}(x)  \quad \quad \: \:  \text{ if }  m< n,
        \end{cases}
\end{align} 
for all integers $m\ge n$, and where $g_u(\cdot) = g(u,\cdot)$.

We say  $g$ has the uniform attraction property (UAP) if for each $\bar{u} \in \overline{U}$ the process $\phi_{\bar{u}}$	is such that 
there exists a sequence of singleton subsets  $\{A_k(\bar{u})\}$ of $X$ so that $\phi_{\bar{u}}(k+1,k,A_k(\bar{u})) = A_{k+1}(\bar{u})$  for all $k \in \mathbb{Z}$ and 
\begin{equation} \label{Ieq_UAP}
	\lim_{j\to \infty} \sup_k \sup_{\bar{u}} dist\Big( \phi_{\bar{u}}(k,k-j,X), A_k(\bar{u})\Big) = 0, 
	\end{equation}
where $dist$ denotes the Hausdorff semi-distance.  More details  behind the definition is available in Section~\ref{Sec_Uniform}. In Theorem~\ref{Thm1} we show that the UAP is equivalent to the USP. When UAP holds it also follows  that $\Psi(\bar{u})$ is a uniform attractor (definition of nonautonomous attractors in \eqref{eq_uni_forward_attractor}) for each $\bar{u}$. Previously, $\Psi(\bar{u})$ was known to be a pullback attractor (e.g., \cite{kloeden2011nonautonomous,manjunath2014dynamics}). Uniform attractors are stronger versions of  nonautonomous attractors like pullback (forward) attractors that are not necessarily attractive when time is run in the negative  (positive) direction (for, e.g., \cite{kloeden2012limitations}).   We also remark that it turns out the  UAP that we have defined is equivalent to having a uniform attractor with fibres that are singleton subsets in the skew-product formalism (see Remark~\ref{Remark_Skew} and Remark~\ref{Rem_UAP}). Under some contraction conditions, existence of uniform attractors in the skew-product formalism can be found in \cite{kloeden2004uniform} and also what are called uniform attracting solutions in \cite{ceni2020echo}. There are other contractive conditions that guarantee USP like the \emph{uniform input forgetting property} \cite[Theorem 6]{grigoryeva2019differentiable} and  the \emph{uniform state contracting property} \cite[Proposition 5]{jaeger2001echo}) but without reference to non-autonomous attractors.

The UAP is not a mere technical point, though. In the scenario where solutions of $g$ are to be computed for information processing, we initialize the driven system $g$ satisfying the UAP with an arbitrary initial value $z_m \in X$, then the sequence $z_{m+1}, z_{m+2}, z_{m+3},...$ satisfying $z_{k+1}= g(u_k,z_k)$  for $k \geq m$ evolves to get within the $\epsilon$ distance of a component of the solution in less than $N(\epsilon)$ time units.   Thus when a solution needs to be computed, by omitting a few values of $z_i$ (referred to as washing out initial conditions in the literature e.g., \cite{jaeger2004harnessing}), for practical purposes the succeeding values of $z_i$ are indistinguishable from the actual solution components. The UAP condition also gives rise to numerical stability in an application involving forecasting \cite{manjunath2021universal}.

\begin{center} 
\begin{figure}[t]
\centering
\includegraphics[width=9cm]{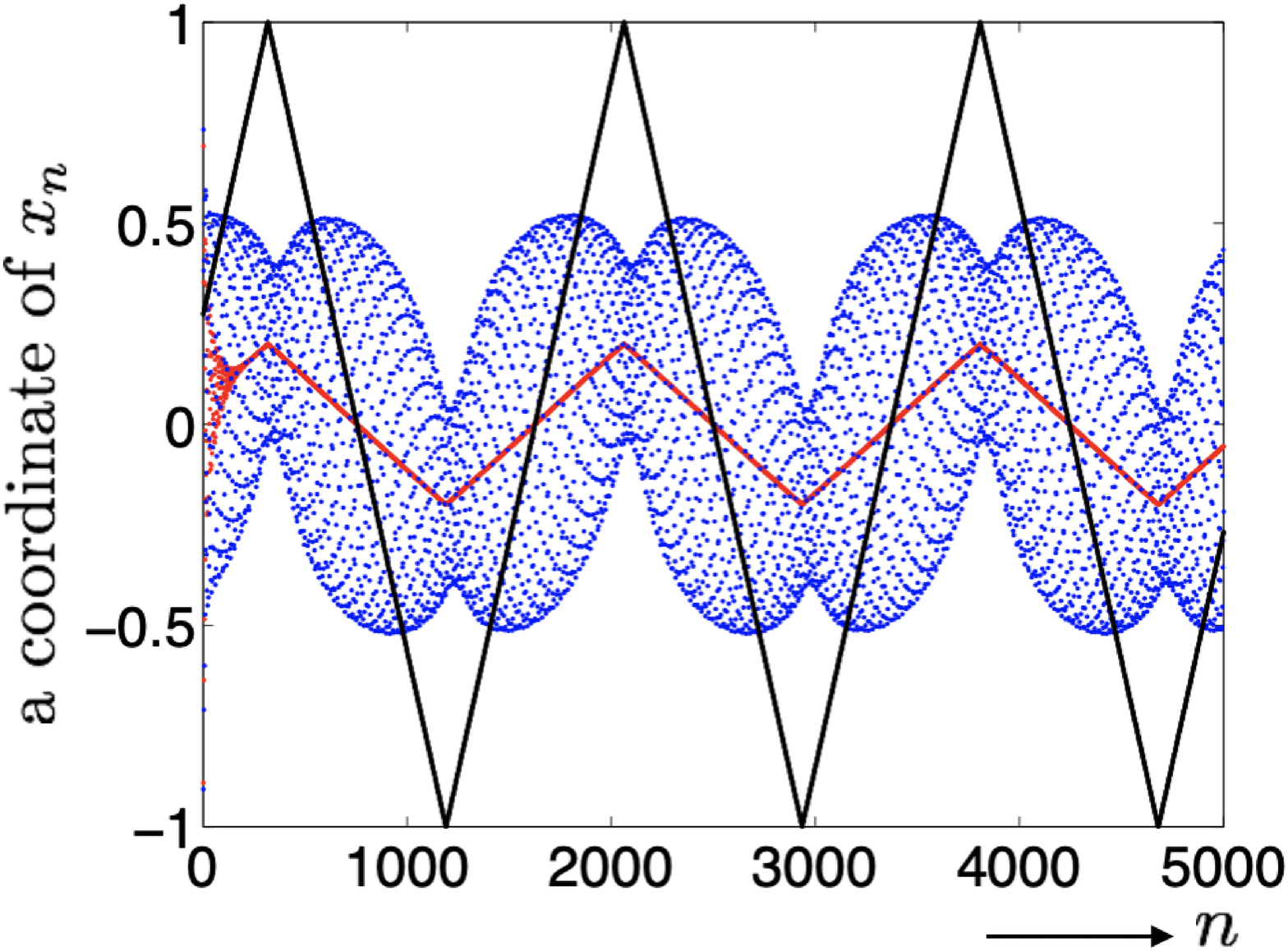}
\caption{A coordinate of a simulated solution $(x_0,x_1,\ldots,x_{5000})$ of a driven system in \eqref{eq_RNN} in Example~\ref{Ex_RNN}  plotted in red (with parameter $\alpha=0.99$ (where $g$ has the USP) and blue ($\alpha=1.05$ where $g$ does not have the USP) while the matrices $A$ and $B$ are randomly generated.}
\label{Fig_Non_Complexity}
 \end{figure}
\end{center}

To explain the importance of the second main result (Theorem~\ref{Thm2}) in the work, we consider the following fundamental question: is the `temporal complexity of the solution of $g$ solely determined by the `complexity of the input'? In other words, would the sequence of mappings $\{g_{u_{n}}(\cdot)\}$ not contribute to the complexity in the solution. 
Continuing to explain the question further, we consider a numerical experiment. An example of a numerical simulation of a solution of a driven system (as in Example~\ref{Ex_RNN}) with two different parameters (see Fig.~\ref{Fig_Non_Complexity}) so that $g$ has the USP for only one of the parameters.  A coordinate of a simulated solution is plotted in red and blue for the two parameters and the input sequence is shown in black. As it may be observed, the coordinate of the solution shown in red (when $g$ has the USP) ``locks" onto a trajectory due to the UAP, and seems to just follow the input in its temporal variation, while that in blue (when $g$ does not have the USP) has wild behaviour with an oscillatory envelope. We observe that the driven system has introduced new additional complexity to the solution indicated in blue that was not there in the input. 

\vspace{-2cm}

So to make the above question formal through the notion of semi-conjugacy, we first consider
 the \textit{reachable set} of the driven system $g$ to be the union of all the elements of all the solutions, i.e., 
$$X_U :=\Big \{x \in X:  x = x_k \mbox{ where  $\{x_n\}$  is a solution for some  $\bar{u}$} \Big \}.$$ 
We then next consider the space of left-infinite sequences  $\overleftarrow{U}$ (notation: $\overleftarrow{Y} := \prod_{i=-\infty}^{-1} Y_i$, where $Y_i \equiv Y$ and is equipped with the product topology). 

\begin{center}
\begin{figure}[t]
\centering
\includegraphics[width=16cm]{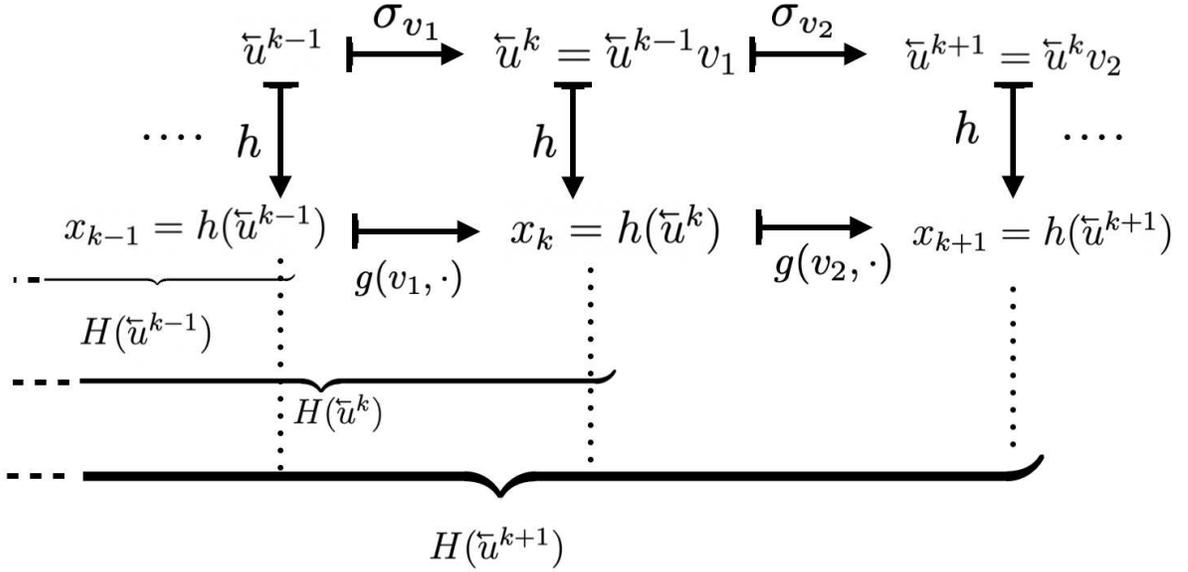}
\caption{Commutating ladder to explain the causal mapping.}
\label{Fig_Causal}
 \end{figure}
\end{center}

We consider $\cev{u}^{n}:=(\ldots,u_{n-2} ,u_{n-1}) \in \overleftarrow{U}$ an input  up to time $n-1$, and define the operation of appending a new input value $v\in U$ at time $n$ by the map $\sigma_v:   
\cev{u}^{n} \mapsto \cev{u}^{n}v$, 
where symbolically  $\cev{u}^{n}v:=(\ldots,u_{n-2},u_{n-1}, v)$, and then ask a fundamental question as to what exactly is required for the dynamics of the driven system to be topologically semi-conjugate to $\sigma_v$ (so that the solution has no new additional new 'complexity)? This is equivalent to asking what exactly $g$ should satisfy for the existence of a continuous surjective map $h: \overleftarrow{U} \to X_U$ (when $X_U$ inherits the subspace topology) so that the following diagram commutes for any $v\in U$: 
\begin{equation}  \label{comm_h}
    \psset{arrows=->, arrowinset=0.25, linewidth=0.6pt, nodesep=3pt, labelsep=2pt, rowsep=0.7cm, colsep = 1.1cm, shortput =tablr}
 \everypsbox{\scriptstyle}
 \begin{psmatrix}
 \overleftarrow{U} & \overleftarrow{U}\\%
 X_U & X_U
 \ncline{1,1}{1,2}^{\sigma_v} \ncline{1,1}{2,1} <{h}
 \ncline{1,2}{2,2} > {h}
 \ncline{2,1}{2,2}^{g(v,\cdot)}.
 \end{psmatrix}
\end{equation} 
We call any such $h$ a universal semi-conjugacy. 
It can be established that the USP implies the existence of a function $h$ 
in a more straightforward way (e.g., \cite{Manju_IEEE}), but in Lemma~\ref{Lemma_imp} we show that the existence of a universal semi-conjugacy is equivalent 
to the USP. We also observe that  $h(\cev{u}^k) = x_k$, where $x_k$ is the value of the solution $\{x_n\}$ at the $k^{\mathit{th}}$ instant for any input $\bar{u}$ whose left-infinite segment is $\cev{u}^k$.  The positive feature of this result is that we can not only show that we can embed a left-infinite input in $\overleftarrow{X}_U$ but also prove the existence of a map $H:  \overleftarrow{U} \to 
\overleftarrow{X}_U$ when ($\overleftarrow{X}_U$ has the product topology) 
have the following diagram commute if and only if $g$ has the USP (Theorem~\ref{Thm2}): 
\begin{equation} \label{comm_H}
\psset{arrows=->, arrowinset=0.25, linewidth=0.6pt, nodesep=3pt, labelsep=2pt, rowsep=0.7cm, colsep = 1.1cm, shortput =tablr}
 \everypsbox{\scriptstyle}
 \begin{psmatrix}
 \overleftarrow{U} & \overleftarrow{U}\\%
 \overleftarrow{X}_U & \overleftarrow{X}_U.
 \ncline{1,1}{1,2}^{\sigma_v} \ncline{1,1}{2,1} <{H}
 \ncline{1,2}{2,2} > {H}
 \ncline{2,1}{2,2}^{\tilde{g}_v}
 \end{psmatrix}
 \end{equation}
Here $\tilde{g}_v: (\cdots,x_{-2},x_{-1}) \mapsto (\cdots,x_{-2},x_{-1},g(v,x_{-1}))$.  We find in our results that $H(\cev{u}^n) =  (\ldots,h(\cev{u}^{n-1}),h(\cev{u}^n)))$.
We call $H$ a causal mapping, and if in addition it also embeds $\overleftarrow{U}$ in $\overleftarrow{X}$ as a causal embedding. Thus when a  causal embedding $H$ exists we have embedded the (left-infinite) input space $\overleftarrow{U}$ onto the (left-infinite) solution space $\overleftarrow{X}_U$ and this is not just an embedding, but the temporal information in the input is inherited in the solution due to the commutativity in \eqref{comm_H}.  A summary of the action of the mappings $h$ and $H$ when $g$ has the USP  is illustrated in Fig.~\ref{Fig_Causal}. During the review of this paper, an interesting question from a referee motivated us to prove (see Lemma~\ref{Lemma_imp} and Proposition~\ref{Prop_Discont}) that if suppose the commutativity in \eqref{comm_h} and \eqref{comm_H} holds without assuming the continuity of $h$ and $H$,  then such $h$ and $H$ have to be discontinuous when $g$ does not have the USP. This gives a plausible explanation to the empirical sensitivity of solutions to a very small amount of noise in the input in driven systems without the USP. Such an experiment can be found in \cite[Figure 3]{Manju_2020}.

We remark out that in systems where time flows continuously rather than discretely as we have considered here, then there are analogs of the functions $h$ and $H$, and are respectively called the \emph{functionals} and \emph{operators} in \cite[Section 2.1]{boyd1985fading}.  In particular, the commutativity in \eqref{comm_H} without reference to the surjectivity is called time-invariance, and when $X=\mathbb{R}$ approximations to the \emph{operator } can be obtained by what is called a Volterra series operator \cite{boyd1985fading}.

In practice, inputs are often restricted to subsets of $\overleftarrow{U}$.  As an application of causal embedding for forecasting dynamical systems \cite{manjunath2021universal}, we point out the finite set of components functions $H_2(\cev{u}^n):= (h(\cev{u}^{n-1}),h(\cev{u}^n))$ can be used to embed such subsets into $X \times X$ (see Remark~\ref{Rem_Partition} for details).

Our results on the driven system do not deal with reconstructing hidden information as in the Takens embedding theorem, but we point out that we can find a $g$ with the USP to implement the delay coordinates $\Phi_{2d,\theta}$ in the Takens embedding theorem. We do that after restating the Takens embedding  theorem from \cite{takens1981detecting}: 
\begin{Theorem*}
	 Let $W$ be a compact manifold of dimension $m$, and $d\ge m$ so that $2d$ is an integer. It is a 
            generic property for the pair $(T, \theta)$,  where $T:W \to W$ is
            a smooth diffeomorphism, and $\theta:W \to \mathbb{R}$ a smooth function, the map $\Phi_{2d,\theta}:W \to \mathbb{R}^{2d+1}$ defined on $W$ by 
            $\Phi_{2d,\theta}(w) := (\theta(T^{-2d}w)\ldots,\theta(T^{-1}w),\theta(w))$
            is a diffeomorphic embedding; by `smooth' we mean at least $C^2$. Consequently, there exists a map $F_\theta: \Phi_{2d,\theta}(W) \to \Phi_{2d,\theta}(W)$ defined by $$F_\theta: (\theta(T^{-2d}w),\ldots,\theta(T^{-1}w),\theta(w)) \mapsto 
            (\theta(T^{-2d+1}w),\ldots,\theta(w),\theta(Tw))$$
           so that $(W,T)$ is topologically conjugate to 
            $(\Phi_{2d,\theta}(W), F_\theta)$. 
\end{Theorem*}
We next recall a ``linear" driven system as in \cite{grigoryeva2021learning} with $U = \theta(W)$ and $X= \Phi_{2d,\theta}(W)$ and $g(u,x):= uC + Ax$, where $A$ is a lower shift matrix (a binary matrix with ones only on the subdiagonal) of dimension $2d+1$, and $C$ is the transpose of the vector $[1,0,\cdots,0]$. Then after initializing the system with any initial condition $x_0$, the system's states are identical to that of the delay-coordinates $\Phi_{2d,\theta}(w)$ in not more than $2d+1$ time-steps. Such a system has the USP since $\Psi(\{\theta(w_i)\}) = \{\Phi_{2d,\theta}(w_i)\}$, where $w_{i+1}=T(w_i)$.

In the final main result in Theorem~\ref{Thm3}, we show that in the case of $g$ having the USP, the solution-map $\Psi$ is a continuous mapping (embedding) whenever $H$ is a causal mapping (embedding).  We go beyond the study of the continuity of a single-valued mapping like $h$ or $H$ and study how a collection of solutions behave as a function of the input as in \cite{Manju_2020}. Towards that end, we consider the reachable states at a \emph{ given time} as a set-valued function of $\cev{u}$. With the space of subsets of such reachable states being equipped with the Hausdorff metric, we describe the continuity of such a set-valued function \emph{at a given time} as a notion of stability which we call the \emph{global-input related stability property} (GIRST) (see Definition~\ref{Def_Girst}) and relate it to the USP. In particular, we show (Proposition~\ref{Prop_GIRST}) that the USP implies GIRST, and conversely, if $g$ has the GIRST we show that if for some input $\bar{u}$ there exists only one solution, then $g$ should have the USP.

We finally point out that the results in this paper can also be interpreted in a study of the perturbation of an autonomous dynamical system. Suppose $U$ were to contain a single point,   $g$ as in \eqref{eqn1} represents an autonomous dynamical system.  The basin of a stable fixed point of an autonomous dynamical system under small perturbations often persists as an attractive neighborhood of a stable or attracting set, although the stable fixed point may not exist any longer in the perturbed system. In this case, we can model the resultant dynamics by expanding $U$ to be a set with more than one value in \eqref{eqn1} to allow the perturbation to be sequences in $U$.   Without the description of attractors in nonautonomous systems (e.g., \cite{kloeden2004uniform}), the folklore understanding remains that when the perturbation is small in amplitude, the dynamics remains close to the fixed point that existed in the absence of the perturbation. The same lack of understanding holds when a local irreducible attractor (precisely, an attractor contained in the chain-recurrent set (e.g., \cite{akin2010general}) of an autonomous system is perturbed.
   Since the dynamics in the neighborhood of an attracting set is contracting, the USP is a natural assumption on the perturbed system. With the USP assumption, we infer from Theorem~\ref{Thm1} that the nonautonomous trajectory in the perturbed system is a uniform attractor, while Theorem~\ref{Thm2} and Theorem~\ref{Thm3} imply that the complexity of the non-autonomous dynamics cannot be more than that in the perturbation and also there is a one-to-one mapping between the input and the nonautonomous dynamics.

The paper is organized as follows. In Section~\ref{Sec_Prelim}, we present the mathematical definitions and some basic results. In Section~\ref{Sec_Uniform}, we show that the UAP and the USP are equivalent. In Section~\ref{Sec_Embedding}, we present the results connecting the USP to the notions of causal mapping and causal embedding. In Section~\ref{Sec_Stability}, we discuss the global input-related stability notion. In Section~\ref{Sec_Conclusions}, we present a short summary of the results.

\vspace{-3cm}
 
 \section{Preliminaries} \label{Sec_Prelim}
 
{\bf Setting.} A driven system would comprise an input metric space $(U,d_U)$, a compact metric (state) space $(X,d)$ and a continuous function $g: U \times X \to X$. For brevity, we refer to $g$ as a driven system with all underlying entities quietly understood. When $U$ is compact, we call $g$ a compactly driven system.

{\bf Notation.} We denote the
collection of all nonempty closed subsets of the metric space $X$ by $\mathsf{H}_X$. 
When $A$ and $B$ belong to  $\mathsf{H}_X$, the quantity $dist(A,B):= \inf\{ \epsilon : A
\subset B_\epsilon(B)\}$ is the Hausdorff semi-distance between $A$ and $B$,  where
$B_\epsilon(B): = \{ x \in X : d(x,B) < \epsilon\}$ is the open
$\epsilon$-neighborhood of $B$. On $\mathsf{H}_X$ we employ the Hausdorff metric defined by
$d_H(A,B) := \max (dist(A,B) , dist(B,A)) : = \inf\{ \epsilon : A
\subset B_\epsilon(B) \: \: \& \: \:B\subset B_\epsilon(A) \}$.  It is well known that whenever $X$
is a compact metric space, $\mathsf{H}_X$ is also a compact
metric space.   The subspace of $\mathsf{H}_X$ that contains the singleton subsets of $X$ is denoted by $\mathsf{S}_X$. 
We also define the mapping $\mathit{i}: X \to \mathsf{S}_X$ by $ \mathit{i}: x\mapsto \{x\}$.  Slightly abusing notations, when $f$ is any function defined on a space $Y$, then for an $A \subset Y$, we denote $f(A)=\bigcup_{x\in A} f(x)$.

All subspaces of a topological space are equipped with the subspace topology.  If $Y$ is a topological space then we denote the spaces of countable products
$\overline{Y} := \prod_{i=-\infty}^\infty Z_i$ and 
$\overleftarrow{Y} := \prod_{i=-\infty}^{-1} Z_i$
 where $Z_i \equiv Y$ and always equip these spaces with the product topology. Given a bi-infinite sequence $\bar{u}$, we define  $\cev{u}^{n}:=(\ldots,u_{n-2} ,u_{n-1})$ for each $n\in \mathbb{Z}$. Also we represent the concatenation of $\cev{u}^n$ and $v$ by $\cev{u}^{n}v:=(\ldots,u_{n-2} ,u_{n-1}, v)$.
 
 Maps on the product spaces: The map $r:\overleftarrow{U} \to \overleftarrow{U}$ denotes the right-shift map, i.e., $r: (\cdots, y_{-2},y_{-1}) \mapsto 
  (\cdots, y_{-3},y_{-2})$. Given $v\in U$ define $\sigma_v: \overleftarrow{U} \to \overleftarrow{U}$ by 
  $\sigma_v: (\cdots, u_{-2},u_{-1}) \mapsto  (\cdots, u_{-2},u_{-1},v)$. 
Given $v\in U$, the map $\tilde{g}_v: \overleftarrow{X} \to \overleftarrow{X}$ is defined by
$\tilde{g}_v: (\cdots, x_{-2},x_{-1}) \mapsto (\cdots, x_{-2},x_{-1},g_v(x_{-1}))$.

{\bf Continuity of set-valued functions.} Let $X$ be a topological space and let $P_X$ denote the power set of $X$.  Let $Z$ be another topological space. A function $f: Z \to P_X$ is said to be upper semicontinuous at $z$ if for every open set $V$ in $X$ containing $f(z)$ there exists an open neighborhood $W$ of $z$ such that $f(W) \subset V$. 
A function $f: Z \to P_X$ is said to be lower semicontinuous at $z$ if for every open set $V$ in $X$ such that $f(z)\cap V \not= \emptyset$, then the set $\{z: f(z) \cap V \not= \emptyset\}$ is an open neighborhood of $x$.  A map $f$ that is upper semicontinuous for all $z\in Z$ is called upper semicontinuous.  A map $f$ that is lower-semicontinuous for all $z\in Z$ is called lower semicontinuous.  A map that is simultaneously both upper semicontinuous and lower semicontinuous at a point $z\in Z$ is said to be continuous at $z$. We use the fact (e.g.\cite{berge1997topological}) that if a map is  upper semicontinuous at $z$ and is single-valued, i.e., $f(z)$ is single-valued, i.e., a singleton subset of $X$, then $f$ is continuous at $z$.  While $X$ is a compact metric space, a set-valued function $f$ can instead be treated as a regular function taking values in the space of nonempty compact subsets $\mathsf{H}_X$. In this case,
the continuity of $f$ defined using upper semicontinuity and lower semicontinuity above coincides with the continuity of the $\mathsf{H}_X$-valued function $f$, with of course $\mathsf{H}_X$ being equipped with the Hausdorff metric (e.g., \cite[Theorem 1, p. 126]{berge1997topological}).
Given a sequence of subsets $\{A_n\}$ of a metric space $X$ the collection of those $x$ for which each open ball centered at $x$ intersects $A_n$ eventually is called the Kuratowski limit-inferior (the set described on the right-hand side of described in \eqref{eq_Kur}). The collection of those $x$ for which each open ball centered at $x$ intersects infinitely many $A_n$ is called the Kuratowski limit-superior. 
The Kuratowski-limit of the sequence $\{A_n\}$ is said to exist when the limit-superior is identical \cite{kuratowski2014topology} and equal to the limit-inferior  -- so
Kuratowski-limit $\mathop{\mathrm{Lim}}_{n \to \infty} A_{n}$ while it exists can
be thus written as
\begin{align} \label{eq_Kur}
     \mathop{\mathrm{Lim}}_{n \to \infty} A_{n} &= \left\{ x \in X \;\left|\; \limsup_{n \to \infty} d(x, A_{n}) = 0 \right. \right\}, \nonumber \\
         &= \left\{ x \in X \;\left|\; \begin{matrix} \mbox{for all open neighbourhoods } U \mbox{ of } x, \\ U \cap A_{n} \neq \emptyset
 \mbox{ for large enough } n \end{matrix} \right. \right\}. 
\end{align}
The fact that we would make use of in this paper is that when $X$ is compact, and each $A_n \subset X$ is compact, then $\mathop{\mathrm{Lim}}_{n \to \infty} A_{n}$ exists if and only if the limit of $\{A_n\}$ in the Hausdorff metric (e.g., \cite{berge1997topological,kuratowski2014topology}) exists, and both these limits are then identical.

{\bf Process and Encoding.} Given a driven system $g$, and $u\in U$, we denote $g_{u}(x) := g(u,x)$. Clearly $g_u: X \to X$.  Suppose a driven system $g$ has been ``fed" input values $u_{m},u_{m+1},\ldots,u_{n-1}$ in that order starting at time $m$ on all initial values in $X$. Then the map $g$ ``transports" a state-value $x\in X$ at time $m$ to a state-value at time $n$ which is $g_{u_{n-1}}\circ \cdots \circ g_{u_{m}}(x)$.  
Formally, for every choice of  $\bar{u}$ we can define for all $m\le n$, the function that maps a state at $x$ at time $m$ through the inputs $u_{m}, u_{m+1}, \dots u_{n-1}$ to the state at time $n$ given by the composition-operator called a process by several authors (e.g.,\cite{kloeden2011nonautonomous}) by the map $\phi_{\bar{u}}: \mathbb{Z}^2_{\ge} \times X \to X$, where $\mathbb{Z}^2_{\ge} := \{ (n,m) : n \ge m, n,m \in \mathbb{Z} \}$ and
\begin{align*}
  \phi_{\bar{u}}(n,m,x) := \begin{cases}
        x \quad \quad \quad \quad \quad \quad \quad \text\quad \quad \quad \quad \quad \quad \:\text{ if $n=m$,}
        \\
        g_{u_{n-1}}\circ \cdots \circ g_{u_{m+1}} \circ g_{u_{m}}(x)  \quad \quad \: \: \text{ if }  m< n.
        \end{cases}
\end{align*}
Note that the set inclusion $\phi_{\bar{u}}(m+2,m,X) \subset \phi_{\bar{u}}(m+2,m+1,X)$ holds for all $m$ since  $g_{u_{m+1}} \circ g_{u_{m}}(X) \subset g_{u_{m+1}}(X)$.  Based on this observation, it follows that  $\phi_{\bar{u}}(n,m,X)$ is a decreasing sequence of sets, i.e., $\phi_{\bar{u}}(n,m-1,X) \subset \phi_{\bar{u}}(n,m,X)$ for an $m$. Hence, if the entire left-infinite input $\{u_m\}_{m<n}$ had influenced the dynamics of the driven system $g$, the system would have evolved at time $n$ to one of the states in the nested intersection
\begin{equation} \label{Seqn_Xn}
X_n(\bar{u}) := \bigcap_{m<n} \phi_{\bar{u}}(n,m,X).
\end{equation}
When $X$ is compact,  $\phi_{\bar{u}}(n,m,X)$ is a closed subset of $X$, and hence $X_n(\bar{u})$ being a nested-intersection is a nonempty closed subset of $X$. Thus $X_n(\bar{u})$ belongs to $\mathsf{H}_X$ and is also the Kuratowski limit (see \eqref{eq_Kur}) i.e.,
$
X_n(\bar{u}) = \{x \in X : \lim_{m\to \infty}d(x,\phi_{\bar{u}}(n,n-m,X))=0 \}
$. Further, the sets $\phi_{\bar{u}}(n,n-m,X)$ are all compact in the compact space $X$, this limit coincides in the Haussdorff metric, i.e., 
\begin{equation} \label{Seqn_Xn2}
X_n(\bar{u}) = d_H\,\mhyphen\lim_{m\to \infty}  \phi_{\bar{u}}(n,n-m,X).
\end{equation}
The set $X_n(\bar{u})$ being the set of all reachable states at time $n$ has another interesting feature, it is the union of the $n^{\text{th}}$ coordinate projection of all solutions when the input to the driven system $g$ is $\bar{u}$. A proof of this simple result can be found in the literature (e.g., \cite[Lemma 2.15]{kloeden2011nonautonomous} or \cite[Lemma 2.2]{Manju_RDE}) or \cite[Proposition 5.5]{kloeden2013discrete}). The family of sets $\{X_n(\bar{u})\}$ is called the representation of $\bar{u}$ in \cite{Manju_2020}

\begin{Lemma} \rm \label{Lemma_soln}
	Let $g$ be a driven system. Then $x\in X_k(\bar{u})$ if and only if there is a solution $\{x_n\}$ with the input $\bar{u}$ such that $x_k=x$.
\end{Lemma}

\begin{Remark} \label{Rem_Contraction} \rm 
 In view of the definition of USP in Section~\ref{Sec_Intro} and Lemma~\ref{Lemma_soln},  $g$ has the USP if and only if 
$X_n(\bar{u})$ is a singleton subset of $X$ for all $\bar{u} \in \overline{U}$ and $n \in \mathbb{Z}$ -- the nested intersection in \eqref{Seqn_Xn} topologically contracts to a single point. Hence, succinctly, we can  say $g$ has the USP if and only if  $X_0(\bar{u}) \in \mathsf{S}_X$ for all $\bar{u} \in \overline{U}$.  This is since if $X_n(\bar{v})$ is not a singleton subset for some $\bar{v}$, then by defining $u_k = v_{k-n}$, then  $X_0(\bar{u})$ is also not a singleton subset.
\end{Remark}

Note that by definition, $X_n(\bar{u})$ is a function of 
$\cev{u}^n$ alone i.e., $(u_n,u_{n+1},\ldots)$ would not influence $X_n(\bar{u})$.  In particular, for $n=0$, i.e., $\bigcap_{m<0} \phi_{\bar{u}}(0,m,X)$ is a function of $\cev{u}$ alone. Henceforth, for convenience, we denote $X_0(\bar{u})$ and its ``finite-time" approximation  (in view of \eqref{Seqn_Xn2}) by $\mathcal{E}(\cev{u})$ and $\mathcal{E}_n(\cev{u}) = \phi_{\bar{u}}(0,-n,X)$ respectively. 
Without further mention, whenever $\bar{u}$ appears concurrently with $\cev{u}$, $\bar{u}$ denotes some fixed right-infinite extension of $\cev{u}$ so that the time-indices are preserved, i.e., $\cev{u}$ contains elements with all negative indices beginning from $-1$. For the convenience of the reader, we tabulate these as set-valued functions of the input:  \\
\renewcommand{\arraystretch}{1.4}
\begin{tabular}{|c|l|l|}
\hline 	 Map  & Definition & Remarks\\  \hline 
 $X_n: \overline{U} \to \mathsf{H}_X$  & $X_n(\bar{u}) = \bigcap_{m<n} \phi_{\bar{u}}(n,m,X)$ & set of reachable states at time $n$\\ \hline
  $\mathcal{E} : \overleftarrow{U} \to \mathsf{H}_X$ &  $\mathcal{E}(\cev{u}) = X_0(\bar{u})$  & set of reachable states at time $0$ \\ \hline
 $\mathcal{E}_n : \overleftarrow{U} \to \mathsf{H}_X$ &  $\mathcal{E}_n(\cev{u}) := \phi_{\bar{u}}(0,-n,X)$ & \pbox{20cm}{set of all reachable states at time $0$ \\ having evolved from time $-n$ \\ i.e., a finite time approximation of $\mathcal{E}$.} \\ \hline
\end{tabular}

We call the map $\mathcal{E}$ the \emph{encoding} map, and in view of Remark~\ref{Rem_Contraction}, we note that $g$ has the USP if and only if $\mathcal{E}(\cev{u}) \in \mathsf{S}_X$ for all $\cev{u} \in \overleftarrow{U}$. This equivalence helps us simplify the arguments in the proofs. We consider the following examples:
\begin{Example} \rm  \label{Example_1} 
	We borrow an example from \cite{kloeden2012limitations}. Let $g(u,x) =\frac{ux}{1+|x|}$ with $U =[0,2]$ and $X=[-2,2]$. Clearly, $g$ is a driven system since $g: U \times X \to X$. Consider $\cev{u} := (\ldots,2,2)$. 
	Then $\mathcal{E}(\cev{u}) = [-1,+1]$ since $g_2([-1,1]) = g(2,[-1,1]) = [-1,1]$ 
	since $-1$ and $+1$  are  fixed points and iterates of $g_2(x)$ for $x\notin [-1,1]$ converge to one of these points.  Since $\mathcal{E}(\cev{u}) \notin \mathsf{S}_X$, $g$ does not have the USP. The reader may note that in view of Lemma~\ref{Lemma_soln}, there is a solution passing through each point $[-1,1]$ at any given point in time as well. 
		\end{Example}

\begin{Example} \label{Ex_RNN} \rm
Let	$U \subset \mathbb{R}^K$ be any compact subset 
and $X=[-1,1]^N$ (the product of $N$ copies of $[-1,1]$) and $g$ be defined by
\begin{equation} \label{eq_RNN}
	g(u,x) = \overline{\tanh}(Au + \alpha Bx),
\end{equation}
where $A$ and $B$ are real matrices with dimension $N \times K$ and $N \times N$ respectively, $\alpha>0$ is a real parameter, and $\overline{\tanh}(*)$ is (the nonlinear activation) $\tanh$ performed component-wise on $*$. It is proved in \cite[Theorem 2]{Manju_ESP} that
when the spectral norm of $\alpha B<1$, then $g$ has the USP.  
\end{Example}

We also note that in general, for any driven system $g$  the encoding map is upper-semicontinuous. We make use of this fact later, but nevertheless, we provide proof of it here for completeness. 

\begin{Lemma} \label{Lemma_USC2}
Let $g$ be a driven system.  Then
the input-encoding map $\mathcal{E}: \overleftarrow{U} \to \mathsf{H}_X$ is upper-semicontinuous. 
\end{Lemma}

{\bf Proof. } Fix $\cev{u}$. Let $\mathcal{E}(\bar{u}) \subset V$ where $V$ is open in $X$.   Fix any choice of $\bar{u}$ so that its left-infinite sequence matches with $\cev{u}$. 

Since $\mathcal{E}(\bar{u})  = d_H\,\mhyphen\lim_{j\to \infty}  \phi_{\bar{u}}(0,-j,X)$, by the definition of the limit, there exists an $n$ such that  $\phi_{\bar{u}}(0,-n,X) \subset V$. Fix such an $n$. To prove the lemma it is sufficient to show that there exists a neigbhorhood $W$ of $\bar{u}$ such that $\phi_{\bar{v}}(0,-n,X) \subset V$ for all $\bar{v}\in W$ since 
 we know that the sequence obtained by increasing $n$ in $\phi_{\bar{u}}(0,-n,X)$ is  decreasing in $X$. Denote $U^{(-n,-1)}$ to be the product space $\prod_{i=-1}^{-n} Z_i$, where $Z_i$ is $U$.

Define $f(\mathbf{u},\cdot):= \phi_{\bar{u}}(0,-n,\cdot)$ where $\mathbf{u} := (u_{-n},\ldots, u_{-1})$. Hence $f: U^{(-n,-1)} \times X \to X$. Since $\phi_{\bar{u}}(0,-n,\cdot)$ is continuous as it is a composition of continuous functions due to its definition in \eqref{eq_Process}, and $g$ is continuous in $u$, it follows that $f$ is continuous.  Since $\phi_{\bar{u}}(0,-n,X) \subset V$ we have $f(\mathbf{u},x) \subset V$ for all $x\in X$.  Now since $f$ is continuous for each $x\in X$, there exists an open set $E_x \subset X$ containing $x$ and $W_x \subset U^{(-n,-1)}$ containing $\mathbf{u}$ such that $f(W_x,E_x) \subset V$. The collection $\{E_x\}$ forms an open cover of $X$, and since $X$ is compact we can find an open cover $\{E_{x_{1}}, \ldots, E_{x_{K}}\}$ such that $f(W_{x_{i}}) \subset V$ for all $1 \le i \le K$.  The finite intersection $W^n = \bigcap_{i=1}^K W_{x_{i}}$ is open in 
$U^{(-n,-1)}$. Hence $f(W^n,E_{x_{i}}) \subset V$ for all for all $1 \le i \le K$ and since $\{E_{x_{1}}, \ldots, E_{x_{K}}\}$ covers $X$, $f(W^n,X) \subset V$. 

Let $W = \{\bar{v} \in \overline{U} : (v_{-n},\ldots, v_{-1}) \subset W^n \}$. Clearly, $W$ is open in $\overline
{U}$ and 
since for $\bar{v}\in W$, $\phi_{\bar{v}}(0,-n,x)= f((v_{-n},\ldots, v_{-1}),x)$, we have  $\phi_{\bar{v}}(0,-n,X) \subset V$ owing to the inclusion $f(W^n,E_{x_{i}}) \subset V$.  $\blacksquare$

\begin{Remark} \rm \label{Rem_cont}
When $\mathcal{E}(\cev{u})$ is a singleton subset of $X$, i.e., when $\mathcal{E}(\cev{u}) \in S_X$, it follows easily (see \cite[Theorem 3.1]{Manju_2020}) that the encoding $\mathcal{E}$ is lower semi-continuous as well at $\cev{u}$, and hence $\mathcal{E}(\cev{u})$ is continuous at $\cev{u}$. Thus when the $g$ has the USP,  $\mathcal{E}$ is continuous. 
\end{Remark}

\begin{Remark} \rm \label{Remark_Skew}  Some of the results in the paper can also be proved using the skew-product formalism (e.g., \cite{kloeden2011nonautonomous,kloeden2012limitations}) that comprises an autonomous dynamical system defined by a dynamical system on a base space, and a cocycle mapping.  To implement such a formalism in our case, one can consider the base space to be $\overline{U}$ and the dynamical system defined by the shift map $\sigma: \overline{U} \to \overline{U}$ as  $\sigma(\bar{u}): \{u_n\} \mapsto \{u_{n+1}\}$. By denoting $\pi_0: \bar{u} \mapsto u_0$ and the $n$-fold composition of $\sigma$ by $\sigma^{n}$, the cocycle mapping can be defined on $\mathbb{N}_0 \times \overline{U} \times X \to X$ as
\begin{align} \label{eq_Skew}
  \varphi(n,\bar{u},x) := \begin{cases}
        x \quad \quad \quad \: \: \; \quad \quad \quad \quad \quad \quad \quad \text\quad \quad \quad \quad \quad \quad \:\text{ if $n=0$,}
        \\
         g_{\pi_{0}(\sigma^{n-1}(\bar{u}))}\circ \cdots \circ g_{\pi_{0}(\sigma(\bar{u}))} \circ g_{\pi_{0}(\bar{u})}(x)    \quad \: \:  \: \: \text{ if }  n > 0.
        \end{cases}
\end{align} 
Specifically  $\{\varphi(n,\cdot,*)\}_{n\in \mathbb{N}_0}$ forms a semi-group of continuous functions since $g$ is continuous, and satisfies the so-called cocycle property: $\varphi(n+k,\bar{u},x) = 
\varphi(n,\sigma^{k}(\bar{u}),x)$ for all $k,n \in \mathbb{N}_0$, $\bar{u} \in \overline{U}$ and $x\in X$.
\end{Remark}

\subsection{Nonautonomous invariant sets and attractors} 

Each input $\bar{u}$ induces a nonautonomous dynamical system which is a sequence of self-maps $\{g_{u_{n}}\}_{n\in \mathbb{Z}}$, the dynamics of which is obtained by the update equation $x_{n+1} =g_{u_{n}}(x_n)$ or equivalently in the process notation  by $x_{n+1}=\phi_{\bar{u}}(n+1,n,x_n)$. 
To have a meaningful notion of invariant sets and attractors for nonautonomous systems, they are considered to be a particular collection of time-indexed subsets; see e.g., \cite{kloeden2011nonautonomous,kloeden2013discrete} for a thorough discussion.  Since the state space is compact in our case, we adopt the definitions of different types of nonautonomous attractors in the framework of a process without reference to bounded sets.

\begin{Definition} \rm
Given a driven system, let $\phi_{\bar{u}}$ be the process that is obtained for an input $\bar{u}$ (see \eqref{eq_Process}). A sequence of subsets $\mathcal{A} = \{A_n(\bar{u}) \}_{n\in \mathbb{Z}} \subset X$ is said to be
a $\phi_{\bar{u}}$-invariant set if it satisfies 
$$\phi_{\bar{u}}(n+1,n,A_n(\bar{u})) =A_{n+1}(\bar{u}) \mbox{ for all } n \in \mathbb{Z}.$$\end{Definition}

Nonempty invariant sets always exist when $X$ is compact for a nonautonomous system \cite{kloeden2011nonautonomous}. For driven systems, it can be shown that (see Lemma~\ref{Lemma_inv}) for any process $\phi_{\bar{u}}$, the representation $\{X_n(\bar{u})\}$ of $\bar{u}$ is an $\phi_{\bar{u}}$-invariant set.

 \vspace{-0.5cm}
Nonautonomous attractors for a process obtained by a driven system and an input $\bar{u}$ are a special type of $\phi_{\bar{u}}$-invariant sets. There are different types of nonautonomous attractors for processes depending upon the attractivity from the distant past (pullback attractors) or attractivity into the distant future (forward attractors) is considered, and they are not equivalent in general \cite{kloeden2011nonautonomous}. The attractivity is described by the closeness of sets as measured by the Hausdorff semi-metric $\mbox{dist}$.

\begin{Definition} \rm \label{Def_local}
Given a driven system, let $\phi_{\bar{u}}$ be the process that is obtained for an input $\bar{u}$ and $\mathcal{A} = \{A_n(\bar{u}) \}$ be
a $\phi_{\bar{u}}$-invariant set such that each $A_n$ is compact and $\subset X_n(\bar{u})$. If the
following conditions
\begin{gather}
\lim_{j \to \infty} \mbox{dist}(\phi_{\bar{u}}(n,n-j, X),A_{n}(\bar{u}))   =  0 \mbox { for all } n, \label{eq_attractor_def1}\\
\lim_{j \to \infty} \mbox{dist}(\phi_{\bar{u}}(n+j,n, X), A_{n+j}(\bar{u}))  =  0 \mbox { for all } n, \\
\lim_{j \to \infty} \sup_n  \mbox{dist}(\phi_{\bar{u}}(n,n-j, X), A_n(\bar{u}))   =   0,\label{eq_uni_pull_attractor} \\
\lim_{j \to \infty} \sup_n \mbox{dist}(\phi_{\bar{u}}(n+j,n, X), A_{n+j}(\bar{u}))  =  0 \label{eq_uni_forward_attractor},
\end{gather}
holds, then in that order, $\mathcal{A}(\bar{u})$ is respectively called a
pullback-attractor, forward-attractor, uniform-pullback-attractor and uniform-forward-attractor of the input $\bar{u}$ or of the process $\phi_{\bar{u}}$. In  particular, if $\{A_n(\bar{u})\}$ is contained in $\mathsf{S}_X$ then they  also called \emph{point attractors} of their types.
\end{Definition}

The following result is known and we prove it here for completeness. The proofs are adopted from the results in \cite{Manju_RDE} and \cite{kloeden2011nonautonomous}.

\begin{Lemma} \label{Lemma_inv}
Let $g$ be a driven system. Then the representation $\{X_n(\bar{u})\}$ of $\bar{u}$ is an $\phi_{\bar{u}}$-invariant set. Hence $\phi_{\bar{u}}(n,n-j,X_{n-j}(\bar{u})) = X_n(\bar{u})$  for all $j,n \in \mathbb{Z}$.  Also, the representation $\{X_n(\bar{u})\}$ of $\bar{u}$ is a pullback attractor.
\end{Lemma}

{\bf Proof. } First we make use of the elementary fact: let  $A_1,A_2,\ldots,$
 be a collection of nonempty closed subsets of a compact space $X$
such that $A_{i+1} \subset A_i$ and
$f:X \to X$ be a continuous function. Then
\begin{equation} \label{eq_intersection2}
f(A) = \bigcap_{i=1}^\infty f(A_i), \mbox{ where } A:=  \bigcap_{i=1}^\infty A_i.
\end{equation} 
This is true since $f(A) \subset \bigcap_{i=1}^\infty f(A_i)$ trivially. When $y \in f(A_i)$ for all $i$ $\Longrightarrow$ there exists  $x_i \in A_i$  such that $f(x_i) = y$ for all $i$. Any limit point $x'$ of $\{x_i\}$  is contained  in $A$ since $A$ is closed. Hence, $x' \in A$. By the continuity of $f$, $f(x') = y$, and so $\bigcap_{i=1}^\infty f(A_i) \subset f(A)$. Next, we prove the $\phi_{\bar{u}}$-invariance by deducing
$\phi_{\bar{u}}(n+1,n,X_n(\bar{u}))  =  \phi_{\bar{u}}(n+1,n, \bigcap_{m<n} \phi_{\bar{u}}(n,m,X)) 
\stackrel{\mbox{(\ref{eq_intersection2})}}{=}   \bigcap_{m<n} \phi_{\bar{u}}(n+1,n,\phi_{\bar{u}}(n,m,X))   =   \bigcap_{m<n} \phi_{\bar{u}}(n+1,m,X)
 =   \bigcap_{m<n} \phi_{\bar{u}}(n+1,m,X) \cap \phi_{\bar{u}}(n+1,n,X) 
 =  \bigcap_{m<n+1} \phi_{\bar{u}}(n+1,m,X), 
 =  X_{n+1}(\bar{u}).$

To show that $\{X_n(\bar{u})\}$ is a pullback attractor, assume on the contrary, i.e., assume 
$
\limsup_{j \to \infty} \mbox{dist}(\phi_{\bar{u}}(n,n-j, X),X_{n}(\bar{u}))   > 0
$ for some fixed $n$. So there exist sequences $x_{j_{k}} \in \phi_{\bar{u}}(n,n-j_{k}, X)$ and $j_{k} \to \infty$ and $\mbox{dist}(x_{j_{k}}, X_n) > \epsilon$ for all $k \in \mathbb{N}$. Without loss of generality assume $x_{j_{k}} \to x_0$ since there always exists a subsequence that converges. So 
 $\mbox{dist}(x_0, X_n) \ge \epsilon$.  However, since $X_n(\bar{u}) = \bigcap_{j>0} \phi_{\bar{u}}(n,n-j,X))$, it contains the limit $x_0$ as it is a nested intersection in a compact space which implies $\mbox{dist}(x_0, X_n)=0$, and this is a contradiction. $\blacksquare$

Conditions on the existence of forward attractors 
in the process formulation were not fully known until a recent result in \cite{kloeden2016construction}. Next, it also turns out that 
if \eqref{eq_uni_pull_attractor} holds then 
\eqref{eq_uni_forward_attractor} also holds and vice-versa. Hence a uniform-pullback-attractor or a uniform-forward-attractor is referred only as a uniform attractor.

\vspace{-0.5cm}

\section{Uniform Attracting Property of a Driven System}  \label{Sec_Uniform}

Here, we bring out the equivalence between the USP and the uniform attraction property for compactly driven systems. For a driven system, each solution $\Psi(\bar{u})$ is a pullback attractor, but for a compactly driven system, it can be shown that it is a uniform attractor. In fact, we can bring out something in common between the uniform attractors that arise for different $\bar{u}$, for compactly driven systems. We first define this commonality between the uniform (point) attractors:

\begin{Definition} \rm \label{Def_UAP}
Let $g$ be a driven system. Then $g$ is said to have the uniform attraction property (UAP) if for each $\bar{u} \in \overleftarrow{U}$ the process $\phi_{\bar{u}}$	has a uniform point attractor $\{A_k(\bar{u})\}$ and 
\begin{equation} \label{eq_UAP}
\lim_{j\to \infty} \sup_k \sup_{\bar{u}} dist\Big( \phi_{\bar{u}}(k,k-j,X), A_k(\bar{u})\Big) = 0.
\end{equation}
\end{Definition}

The uniform attraction property in Definition~\ref{Def_UAP} is obviously stronger than having a uniform point attractor $\{A_k(\bar{u}))\}$ for each $\bar{u}$ since one could obtain an integer $N$ so that every initial condition $x$ evolves to get within the $\epsilon$ distance of a component of the point attractor in less than $N$ time units. 
 In the case of the USP, i.e., when the solution map $\Psi(\bar{u})$ is well-defined, these point attractors are solutions.

\begin{Theorem} \label{Thm1}
	A compactly driven system $g$ has the USP iff $g$ has the UAP. 
\end{Theorem}

{\bf Proof. } (UAP $\Longrightarrow$ USP).  Fix $\cev{u}$. Let $\{A_n(\bar{u})\}$ be a uniform point attractor for $\phi_{\bar{u}}$. Then there is no other solution except $\{\mathit{i}^{-1}(A_n(\bar{u}))\}$ since 
$\{A_n(\bar{u})\}$ is a point attractor. 

(USP $\Longrightarrow$ UAP). If $g$ has the USP,  then $\mathcal{E}(\cev{u}) \in \mathsf{S}_X \subset \mathsf{H}_X$ for all $\cev{u} \in \overleftarrow{U}$.

Consider the sequence $\{\mathcal{E}_n(\cev{u}) \}_{n\ge 1} \subset \mathsf{H}_X$, where $\mathcal{E}_n(\cev{u}) = \phi_{\bar{u}}(0,-n,X)$. Hence, in view of \eqref{Seqn_Xn2},
$\mathcal{E}(\cev{u}) = d_H\text{\,-}\lim_{n\to \infty} \mathcal{E}_{n}(\cev{u})$.

{\bf Claim. }  $\mathcal{E}_n$ converges uniformly to $\mathcal{E}$ on $ \overleftarrow{U}$, i.e., for all $\epsilon > 0$, there exists $N$ such that 
$$
d_H(\mathcal{E}_n(\cev{u}), \mathcal{E}(\cev{u})) < \epsilon
$$
for all $\cev{u} \in \overleftarrow{U}$ and all $n\ge N$. 
Fix $\epsilon >0$. For each $n \in \mathbb{N}$ define the subset of $\overleftarrow{U}$ as:
$$
B_n:= \{ \cev{u} \in \overleftarrow{U}: d_H(\mathcal{E}_n(\cev{u}), \mathcal{E}(\cev{u}))< \epsilon \}.$$
Since $\mathcal{E}_n(\cev{u}) \to \mathcal{E}(\cev{u})$ for each $\cev{u}$, $B_n$ is nonempty for some $n$. The set $B_n$ is open in $\overleftarrow{U}$ for all $n \in \mathbb{N}$ in view of the following facts: (i) $\mathcal{E}_n(\cev{u})$ is decreasing to $\mathcal{E}(\cev{u})$ for all $\cev{u}$; (ii) $\mathcal{E}_n$ is continuous for each $n\in \mathbf{N}$; (iii) $\mathcal{E}$ is continuous owing to taking values in $\mathsf{S}_X$ (iv) $d_H$ being a metric is continuous. Note that $\{B_n\}_{n\in \mathbb{N}}$ is an open cover of $\overleftarrow{U}$, and since $\overleftarrow{U}$ is compact, we have a finite collection of sets $\{B_{n_{1}}, B_{n_{2}}, \ldots, B_{n_{m}}\}$ that covers $\overleftarrow{U}$. Let $N:= \max(n_1,n_2, \ldots , n_m)$. Also since, $\mathcal{E}_n(\cev{u})$ is decreasing, $B_n \subset B_{n+1}$. Therefore $B_N$ covers $\overleftarrow{U}$. 
Thus by definition of $B_N$, we have for all $n \ge N$, $d_H(\mathcal{E}_n(\cev{u}),\mathcal{E}(\cev{u})) < \epsilon$ for all 
$\cev{u} \in \overleftarrow{U}$. Thus the claim is proved. 

Recall that $g$ has the UAP if for all $\epsilon > 0$ there exists $J_\epsilon \in \mathbb{N}$ and $\{A_n(\bar{v}) \} \subset \mathsf{S}_X$ such that 
$$
\sup_k \sup_{\bar{v}} dist( \phi_{\bar{v}}(k,k-j,X), A_k(\bar{v})) \le \epsilon
$$
for all $j \ge J_\epsilon$. Hence to show that $g$ has the UAP, it is sufficient to show given any $\bar{v}$ and any $k \in \mathbb{Z}$, there exists $J_{\epsilon}$ and  $A_k(\bar{v}) \in  \mathsf{S}_X$ such that
$$
dist( \phi_{\bar{v}}(k,k-j,X), A_k(\bar{v})) \le \epsilon
$$ for all $j > J_\epsilon$. Fix $\bar{v}$ and $k\in \mathbb{Z}$ and let $J_\epsilon =N$. Now choose $\bar{u} \in \overline{U}$ such that 
$u_i = v_{k+i}$ for all $i \in \mathbb{Z}$. Hence $(\ldots,u_{-2},u_{-1})=(\ldots, v_{k-2},v_{k-1})$ and thus $\phi_{\bar{v}}(k,k-j,X) = \phi_{\bar{u}}(0,-j,X)$. But 
$\phi_{\bar{u}}(0,-j,X) = \mathcal{E}_j(\cev{u})$. Now let $A_k(\bar{v}) = \mathcal{E}(\cev{u})$. Therefore, 
\begin{eqnarray*}
	dist( \phi_{\bar{v}}(k,k-j,X), A_k(\bar{v})) & \le & d_H (\phi_{\bar{v}}(k,k-j,X), A_k(\bar{v})) \\
	&=&d_H(\mathcal{E}_j(\cev{u}), \mathcal{E}(\cev{u})).
\end{eqnarray*}
But $d_H(\mathcal{E}_j(\cev{u}), \mathcal{E}(\cev{u})) < \epsilon$ for $j > J_\epsilon$ due to the claim above.  Hence the UAP follows. $\blacksquare$

We note that in the skew-product formalism as described in Remark~\ref{Remark_Skew},  the collections of sets $\mathbb{A} = \{A_{\bar{u}} \subset X\}$ where each fibre $A_{\bar{u}}$ is nonempty, compact, and invariant, i.e., $\Psi(n,\bar{u}, A_{\bar{u}}) = A_{\sigma^{n}(\bar{u})}$  is called a uniform attractor (e.g., \cite{kloeden2011nonautonomous}) for the skew-product flow  obtained from the driven system $g$ if 
\begin{equation}  \label{eq_Uni_Skew}
	\lim_{k\to \infty} \sup_{\bar{u}} dist\Big( \varphi(k,\bar{u},X), A_{\sigma^{k}(\bar{u})}\Big) = 0.
\end{equation}

\begin{Remark} \rm \label{Rem_UAP}
Replacing $k$ by $k+j$ in \eqref{eq_UAP} we can restate the UAP condition by \\  $\lim_{j\to \infty} \sup_k \sup_{\bar{u}} dist\big( \phi_{\bar{u}}(k+j,k,X), A_{k+j}(\bar{u})\big) = 0$. Next, by comparing the definitions of a process in \eqref{eq_Process} and the cocycle mapping in \eqref{eq_Skew}, we note that $\phi_{\bar{u}}(k+j,k,x) = \varphi(j,\sigma^k(\bar{u}),x)  
$.  So from \eqref{eq_Uni_Skew}, we 
can conclude that the UAP defined through the process formulation is equivalent to  $\mathbb{A}$ being a uniform attractor in the skew product formalism containing only singleton subsets of $X$ as its fibres.
\end{Remark} 

\vspace{-0.5cm}

\section{Causal Embedding and Conjugacy}

 \label{Sec_Embedding}
Here, we establish the equivalence between the USP and the existence of a causal mapping for a compactly driven system in Theorem~\ref{Thm2}. The key to our results is proving the existence of a universal semi-conjugacy in Lemma~\ref{Lemma_imp}. We also establish conditions under which $H$ is a causal embedding in Theorem~\ref{Thm2}. Another qualitatively different sufficient condition that enables embedding inputs as driven system’s states can be found in \cite[Corollary 23]{grigoryeva2019differentiable}. To the more circumspect reader, we also show in Proposition~\ref{Prop_Discont} any function $H$ that makes \eqref{comm_H} commute, i.e., any $H$ that satisfies Definition~\ref{Def_causal_embedding} except that it need not be continuous would be actually  discontinuous when $g$ does not have the USP.

Recall from Section~\ref{Sec_Intro}, the reachable set for a driven system $g$ is defined as $X_U = \Big \{x \in X:  x = x_k \mbox{ where  $\{x_n\}$  is a solution for some  $\bar{u}$} \Big \}. $.  Also, recall that $X_k(\bar{u})$ is the union of all $k^{\text{th}}$ component of all solutions obtained with the input $\bar{u}$ by Lemma~\ref{Lemma_soln}.
In this light, we can equivalently define  the reachable set $X_U$ by:
\begin{equation}  \label{eqn_reachable}
X_U = \bigcup_{n\in \mathbb{Z},  \:\bar{u} \in \overline{U}} X_n(\bar{u}).
\end{equation}
We also note that for any $x\in X_U$ and $v\in U$, $g(v,x) \in X_U$. This is since if $x\in X_U$, there exists a left-infinite input $\cev{u}^n$ such that $x\in X_n(\bar{u}) = \mathcal{E}(\cev{u}^n)$. So 
by the definition of $\mathcal{E}$ it follows that $x\in\mathcal{E}(\cev{u}^nv)$, and hence $x\in X_U$.   We would use this fact without further mention.

\begin{Definition} \rm \label{Def_causal_embedding}
	Consider a driven system $g$. We say $H$ is a causal mapping if $H: \overleftarrow{U} \to \overleftarrow{X}_U$ is continuous, surjective and satisfies (the commutativity in the diagram in \eqref{comm_H})
	$$
	\tilde{g}_v\circ H = H \circ \sigma_v$$
	for all $v \in \overleftarrow{U}$, where  $\sigma_v : \cev{u}^n \mapsto \cev{u}^nv$ and $\tilde{g}_v: (\cdots,x_{n-2},x_{n-1}) \mapsto (\cdots,x_{n-2},x_{n-1},g(v,x_{n-1}))$.	
	If in addition, $H$ is a homeomorphism, we say $H$ is a causal embedding of $\overleftarrow{U}$ in $\overleftarrow{X}$. 
\end{Definition}

\begin{Definition} \rm
	Consider a driven system $g$. We say $h$ is a universal semi-conjugacy if 
	 $h: \overleftarrow{U} \to X_U$ is continuous, surjective and satisfies (the commutativity in the diagram in \eqref{comm_h}):
	$$
	g_v \circ h = h \circ \sigma_v \mbox{ for all $v\in U$ }.
	$$ 
\end{Definition}

\begin{Lemma} \label{Lemma_agg_limit}
Consider a driven system $g$ and $\phi_{\bar{u}}$ be the process corresponding to an
$\bar{u} \in \overline{U}$. Then,
$$\bigcap_{j=1}^\infty \phi_{\bar{u}}(n,n-j,X) = \bigcap_{j=1}^\infty \phi_{\bar{u}}(n,n-j,X_U).$$
\end{Lemma}

\vspace{-0.5cm}
{\bf Proof. } Since $X_U \subset X$, we have $\phi_{\bar{u}}(n,n-j,X_U)	 \subset \phi_{\bar{u}}(n,n-j,X)$ for all $n,j$, and thus for all $n$,
\begin{equation} \label{eqn_inclusion}
	\bigcap_{j=1}^\infty \phi_{\bar{u}}(n,n-j,X_U)	 \subset \bigcap_{j=1}^\infty  \phi_{\bar{u}}(n,n-j,X).
\end{equation}
Next, by \eqref{eqn_reachable}, we know $X_{n-j}(\bar{u}) \subset X_U$ $\forall j$, so we have $\phi_{\bar{u}}(n,n-j,X_{n-j}(\bar{u})) \subset \phi_{\bar{u}}(n,n-j,X_U)$ for all $j \in \mathbb{N}$. But $\phi_{\bar{u}}(n,n-j,X_{n-j}(\bar{u})) = X_n(\bar{u})$ by Lemma~\ref{Lemma_inv}. 
Hence, $X_n(\bar{u}) \subset \phi_{\bar{u}}(n,n-j,X_U)$ for all $j \in \mathbf{N}$ which implies 
  $X_n(\bar{u}) \subset  \bigcap_{j=1}^\infty\phi_{\bar{u}}(n,n-j,X_U)$.
But $X_n(\bar{u}) =   \bigcap_{k=1}^\infty \phi_{\bar{u}}(n,n-k,X)$ by definition.
   Thus for all $n$,
  	\begin{equation} \label{eqn_inclusion2}
	\bigcap_{k=1}^\infty \phi_{\bar{u}}(n,n-k,X)	 \subset \bigcap_{j=1}^\infty  \phi_{\bar{u}}(n,n-j,X_U).
\end{equation}
By \eqref{eqn_inclusion} and \eqref{eqn_inclusion2}, the lemma follows. $\blacksquare$.

We use the hypothesis of invertibility of $g(\cdot,x): U \to X$ for all $x \in X$ as a sufficient condition for a causal mapping $H$ to be a causal embedding in the statement (ii) of Theorem~\ref{Thm2}. Also 
the reader may note that $g(\cdot,x)$ being invertible for all $x \in X_U$ serves as a sufficient condition as well.

\begin{Example} \rm \label{Example_inv}
It is also easy to verify for the example of the driven system in \eqref{eq_RNN} if $U \subset \mathbb{R}^N$ and $A$ has the same dimension as that of $B$ and is invertible, then $g(\cdot,x)$ is also invertible for all $x \in [-1,1]^N$. 
\end{Example}

\renewcommand{\labelenumi}{(\roman{enumi}).}
\begin{Theorem} \label{Thm2}
	Consider a compactly driven system $g$. 
	\begin{enumerate}
		\item There exists a causal mapping $H: \overleftarrow{U} \to \overleftarrow{X}_U$  if and only if $g$ has the USP;  $X_U$ inherits the subspace topology from $X$. In particular, 
		\begin{equation} \label{eqn_H}
			H(\cev{u}) = \Big(\cdots, h(r^2\cev{u}),h(r\cev{u}),h(\cev{u})\Big),\end{equation} where $h=\mathit{i}^{-1}\circ\mathcal{E}$, where $\mathcal{E}$ is the encoding map, and $\mathit{i}: x \mapsto \{x\}$.

		\item $H$ is a causal embedding when $g(\cdot,x): U \to X$  is invertible for all $x\in X$.
	\end{enumerate}
\end{Theorem}
 
 {\bf Proof.} {\bf Proof of (i).} Let $g$ have the USP, and let $H$ be as in \eqref{eqn_H}. Clearly, $H$ is continuous since both $h: \overleftarrow{U} \to \overleftarrow{X}_U$ by Lemma~\ref{Lemma_imp} and the right-shift map $r: \overleftarrow{U} \to \overleftarrow{U}$ are continuous. Also, $H$ is surjective since $h$ 
again by Lemma~\ref{Lemma_imp} is surjective, and $r$ is also surjective.  
The equality $\tilde{g}_v \circ H = H \circ \sigma_v$ holds if and only if the equality 
 $\pi_j \circ \tilde{g}_v (H(\cev{u})) = \pi_j \circ H(\sigma_v(\cev{u}))$  holds for all $j \in \mathbb{N}$, where $\pi_j$ for the coordinate projections. Since $g$ has the USP, by Lemma~\ref{Lemma_imp}, we note that $h\circ \sigma_v = g_v \circ h$. Then $H :\overleftarrow{U} \to \overleftarrow{X}$ is a causal mapping if the following claim is true.  

{\bf Claim.} $h \circ \sigma_v = g_v \circ h$ $\Longrightarrow$ $H \circ \sigma_v = \tilde{g}_v \circ H$ (which is the required commutativity in the diagram in \eqref{comm_H}). \\
{\bf Proof of Claim.} Suppose $H \circ \sigma_v(\cev{u}) \not= \tilde{g}_v \circ H(\cev{u})$ for some $\cev{u}$. Then $\pi_j(H \circ \sigma_v(\cev{u})) \not= \pi_j(\tilde{g}_v\circ h(\cev{u}))$ for some $j \in \mathbb{N}$. Suppose $j>1$, and since the composition $r^{-j-2}\circ r =r^{-j-1}$,
and since $r \circ \sigma_v(\cev{u}) = \cev{u}$,  $h(r^{-j-2}(\cev{u})) = h(r^{-j-1}\circ \sigma_v(\cev{u})) \not= 
 h(r^{-j-2}(\cev{u}))$ which is absurd. So $j$ has to be equal to $1$. Then $\pi_1(H\circ \sigma_v(\cev{u})) \not= \pi_1(\tilde{g}_v\circ h(\cev{u}))$ means $h\circ \sigma_v(\cev{u}) \not= g_v \circ h(\cev{u})$. The claim is proven. 
 
The converse remains to be proven: that if $H: \overleftarrow{U} \to \overleftarrow{X}$ is a causal mapping, then $g$ has the USP and $H$ will be as in \eqref{eqn_H}. Let 
 \begin{equation} \label{eqn_H2}
 	 H(\cev{u})=(\cdots, \theta_2(\cev{u}),\theta_1(\cev{u})),
 \end{equation}
where $\theta_i: \overleftarrow{U} \to X_U$ is some function that is continuous and surjective. Observe that $\tilde{g}_v \circ H = H \circ \sigma_v$ if and only if $\pi_j \circ \tilde{g}_v \circ H = \pi_j \circ  H \circ \sigma_v$ for all coordinate projections $\pi_j$, $j \in \mathbb{N}$. With $j=1$ and from \eqref{eqn_H2}, $\pi_1 \circ \tilde{g}_v \circ H = g_v \circ \theta_1$ and $\pi_1 \circ  H \circ \sigma_v = \theta_1 \circ \sigma_v$. So we have $g_v \circ \theta_1 = \theta_1 \circ \sigma_v$. So, if for all $\cev{u}$, $g_v \circ \theta_1(\cev{u}) = \theta_1 \circ \sigma_v(\cev{u})$ holds, then $\theta_1$ is a universal semiconjugacy, and thus by Lemma~\ref{Lemma_imp}, 
$g$ has the USP. Further from Lemma~\ref{Lemma_imp}, $\theta_1 = h(\cev{u})$, where $h= \mathit{i}^{-1} \circ \mathcal{E}$. 
This implies $\theta_{k+1} = h(r^k\cev{u})$. This completes the proof of (i). 

{\bf Proof of (ii).} Suppose $H(\cev{u}) = H(\cev{v})$ for some $\cev{u} \not=\cev{v}$. So by \eqref{eqn_H}, $h(r^k\cev{u}) = h(r^k\cev{v})$. 
Since $r^k\cev{u} = (\ldots,u_{-k-2},u_{-k-1})$ for all $k>0$, we have $h(\cev{u}^{-k-1}) = h(\cev{v}^{-k-1})$ for all $k>0$.  Since $\cev{u} \not=\cev{v}$,   fix an integer $m<-1$ such that $u_m\not=v_m$. Also, since $h=\mathit{i}^{-1}\circ \mathcal{E}$, by Lemma~\ref{Lemma_imp}, we have $g(u_k,h(\cev{u}^k))=h(\cev{u}^{k+1})$ for all $k\le-2$. Set $k=m$. Then,
$
g(u_m,h(\cev{u}^m))=h(\cev{u}^{m+1})$ and 
$
g(v_m,h(\cev{v}^m))=h(\cev{v}^{m+1})$ which is not possible since $h(\cev{u}^k) = h(\cev{v}^k)$ for all $k\le 0$ and 
$g(\cdot,x)$ is invertible for all $x\in X$. Therefore $H$ is invertible and a homeomorphism since it's a mapping between two compact Hausdorff spaces. Thus $H:\overleftarrow{U} \to \overleftarrow{X}_U$ is a causal embedding.  $\blacksquare$

\begin{Remark} \rm \label{Rem_Partition}
 In applications of driven systems, the inputs are generally restrictive and they may be drawn only from a subset of the left-infinite sequences  $\overleftarrow{U}$. When inputs are restricted to a subset of $\overleftarrow{U}$, then it may be possible that $H_k(\cev{u}) = \Big(h(r^{(k-1)}\cev{u}),\cdots ,h(r\cev{u}),h(\cev{u})\Big)$ embeds the subset of $\overleftarrow{U}$ into $X_U^{k}$, where $X_U^{k}$ is the space obtained by the product of $k$ copies of $X_U$. An example arises when the input originates from another autonomous dynamical system, i.e., the input is restricted to be a left-infinite orbit from $T:U \to U$, where $T$ is a homeomorphism. In this case,  it is shown in \cite{manjunath2021universal}  the existence of a causal embedding $H$ along with $g(\cdot,x) : U \to X$ being invertible implies that the map $H_2(\cev{u}^n):= (h(\cev{u}^{n-1}),h(\cev{u}^n))$ embeds all left-infinite orbits of $T$ in $X_U \times X_U$, and $H_2$ is also a topologically conjugacy between the single-delay lag dynamics of the driven system and $T$ i.e., there exists a map $G_T:(x_{n-1},x_n) \mapsto (x_{n},x_{n+1})$ that is conjugate to $T$. 
\end{Remark}

 \begin{Lemma} \label{Lemma_imp}
	Consider a compactly driven system $g$.  Then $g$ induces a universal semi-conjugacy $h: \overleftarrow{U} \to X_U$ if and only if $g$ has the USP.  	Further, 
	$h=\mathit{i}^{-1} \circ \mathcal{E}$, where $\mathcal{E}$ is the encoding map, and  $i: x \mapsto \{x\}$. On the other hand suppose  $h: \overleftarrow{U} \to X_U$ is a surjective function so that the diagram in \eqref{comm_h} commutes without being necessarily continuous then such an $h$ is actually not continuous when $g$ does not have the USP. 
	\end{Lemma}

\vspace{-0.5cm}
{\bf Proof. }  Assume $g$ has the USP and then we will show that $h=\mathit{i}^{-1} \circ \mathcal{E}$ is the universal semi-conjugacy. Since $\mathcal{E}$ takes values in $\mathsf{S}_X$, the function $\mathit{i}^{-1} \circ \mathcal{E}$ is well-defined. 
Let  $h(\cev{u}) = \mathit{i}^{-1}(\mathcal{E}(\cev{u}))$. By definition of $\mathcal{E}$, we find (the required commutativity in the diagram in \eqref{comm_h}): 
\begin{eqnarray*}
g_v\circ \mathit{i}^{-1} \circ \mathcal{E}(\cev{u}) & = & \mathit{i}^{-1} \circ \mathcal{E}(\cev{u}v), \\
& =& \mathit{i}^{-1} \circ \mathcal{E}(\sigma_v(\cev{u})).
\end{eqnarray*}
By definition of $X_U$,  $\mathit{i}^{-1} \circ \mathcal{E}: \overleftarrow{U} \to X_U$ is surjective and $\mathcal{E}: \overleftarrow{U} \to \mathsf{H}_X$ is continuous whenever $g$ has the USP (see Remark~\ref{Rem_cont}). The function $\mathit{i}$ is continuous as it is an isometry, and therefore
$\mathit{i}^{-1} \circ \mathcal{E}$ is continuous. Hence, $\mathit{i}^{-1} \circ \mathcal{E}: \overleftarrow{U} \to X_U$ is a universal semi-conjugacy.

We next derive an expression for $\mathcal{E}(\cev{u})$ when $h: \overleftarrow{U} \to X_U$ is surjective, and is not necessarily continuous but satisfies $h \circ \sigma_v = g_v \circ h$ (the commutativity in the diagram in \eqref{comm_h}).  For any $\cev{u} \in \overleftarrow{U}$ fix a $\bar{u}$ whose left-infinite sequence is $\cev{u}$.
We deduce
\begin{eqnarray}
	\mathcal{E}(\cev{u}) &=& d_H \,\mhyphen \lim_{j \to \infty} \phi_{\bar{u}}(0,-j,X), \nonumber \\	&\stackrel{\mbox{(Lemma~\ref{Lemma_agg_limit})}}{=}& d_H \,\mhyphen\lim_{j \to \infty} \phi_{\bar{u}}(0,-j,X_U), \nonumber \\
	&=& d_H \,\mhyphen\lim_{j\to \infty}  \phi_{\bar{u}}(0,-j, h(\overleftarrow{U})), \hspace{0.5cm} (\mbox{since $h$ is surjective}) \nonumber \\
	&=& d_H \,\mhyphen\lim_{j\to \infty} g_{u_{-1}} \circ \cdots \circ g_{u_{-j}} \circ h(\overleftarrow{U}), \hspace{0.5cm} (\mbox{by definition of $\phi_{\bar{u}}$}) \nonumber \\
	&=& d_H \,\mhyphen\lim_{j \to \infty}  g_{u_{-1}} \circ \cdots \circ g_{u_{-j+1}} \circ h(\overleftarrow{U}u_{-j}) \hspace{0.5cm} (\mbox{since } g_{u_{-j}} \circ h(\cev{w}) = h(\cev{w}u_{-j})), \nonumber \\
	&=& d_H \,\mhyphen\lim_{j \to \infty} h(\overleftarrow{U}u_{-j}\cdots 
	u_{-1}), \hspace{0.15cm}(\mbox{applying } g_{u_{-k}}\circ h(\cev{w}) = h(\cev{w}u_{-k}) \mbox{ repeatedly}) \nonumber \\
	&=&d_H \,\mhyphen\lim_{j \to \infty}h(\{(\cev{w}u_{-j}\cdots 
	u_{-1}\}): \cev{w} \in \overleftarrow{U}\}). \label{eqn_alternate_exp} \end{eqnarray}
We note that in the above deduction without assuming the continuity of $h$, we have shown that for each $j$,  the set $\{h(\cev{w}u_{-j}\cdots 
	u_{-1}\}): \cev{w} \in \overleftarrow{U}\}$ is compact since it is equal to $\phi_{\bar{u}}(0,-j,X)$. Since the Kuratowski and the Hausdorff limit coincide for a sequence of compact subsets in a compact space $X$, by  \eqref{eq_Kur}, we use that in \eqref{eqn_alternate_exp} to obtain 
\begin{align} \label{eq_Kur_final}
\mathcal{E}(\cev{u})     
         &= \left\{ x \in X \;\left|\; \begin{matrix} \mbox{for all open neighbourhoods } U \mbox{ of } x, \\ U \cap h(\{(\cev{w}u_{-j}\cdots 
	u_{-1}): \cev{w} \in \overleftarrow{U}\}) \neq \emptyset
 \mbox{ for large enough } n \end{matrix} \right. \right\}.
\end{align}
Suppose if there are two distinct points $x_1$ and $x_2$ of $X_U$ in $\mathcal{E}(\cev{u})$. It follows from \eqref{eq_Kur_final} that there exist two disjoint neighborhoods of $U_1$ and $U_2$ in $X_U$ of $x_1$ and $x_2$ respectively that intersects $h(\{(\cev{w}u_{-j}\cdots u_{-1}: \cev{w} \in \overleftarrow{U}))$ for all $j>J$. Now, this is possible only when $h$ is not continuous at $\cev{u}$. For if, $h$ is continuous, we can find a neighborhood $V$ of $h(\cev{w})$ that is disjoint from at least one of these neighborhoods $U_1$ and $U_2$ and 
\begin{equation} \label{eq_hinclusion} h(\{ \cev{w} u_{-j}\cdots u_{-1} : \cev{w} \in \overleftarrow{U} \}) \subset V \end{equation} for large enough $j$. 
Note that  the set inclusion in \eqref{eq_hinclusion} holds for all large $j$ since 
$\overleftarrow{U}$ is metrizable (e.g., 
$d(\cev{u},\cev{v}) := \displaystyle{\sum_{i=-\infty}^{-1}} d_U(u_i,v_i) / 2^{|i|}$ is a metric) and $\{ \cev{w} u_{-j}\cdots u_{-1} : \cev{w} \in \overleftarrow{U} \}$ can be made to be contained in any open ball of radius $\delta >0$ centered at $\cev{u}$ for all $j$ sufficiently large.  So we have proved if the commutativity in the diagram in \eqref{comm_h} holds, then $\mathcal{E}(\cev{u}) \in \mathsf{S}_X$ (i.e., $g$ has the USP) if $h$ is continuous, and also suppose $\mathcal{E}(\cev{u}) \notin \mathsf{S}_X$ (i.e., suppose $g$ does not have the USP) then $h$ is discontinuous.  
	$\blacksquare$.

\begin{Proposition} \label{Prop_Discont}
Consider a compactly driven system $g$. Suppose there exists a  surjective mapping $H: \overleftarrow{U} \to \overleftarrow{X}_U$ so that \eqref{comm_H} commutes without being necessarily continuous, then $H$ is actually discontinuous when $g$ does not have the USP. 	
\end{Proposition}

{\bf Proof. } Let 
 \begin{equation} \label{eqn_H3}
 	 H(\cev{u})=(\cdots, \theta_2(\cev{u}),\theta_1(\cev{u})),
 \end{equation}
where $\theta_i: \overleftarrow{U} \to X_U$ is some function. Since $H$ is surjective, any $\theta_i$ is also surjective. If the diagram in \eqref{comm_H} commutes, then $\tilde{g}_v \circ H = H \circ \sigma_v$ which 
implies $\pi_1 \circ \tilde{g}_v \circ H = \pi_1 \circ  H \circ \sigma_v$ for the coordinate projection $\pi_1$.  Using \eqref{eqn_H3}, it follows from \eqref{comm_H} that $\pi_1 \circ \tilde{g}_v \circ H = g_v \circ \theta_1$ and $\pi_1 \circ  H \circ \sigma_v = \theta_1 \circ \sigma_v$. So we have $g_v \circ \theta_1 = \theta_1 \circ \sigma_v$. So, if for all $\cev{u}$, $g_v \circ \theta_1(\cev{u}) = \theta_1 \circ \sigma_v(\cev{u})$ holds, then by setting $h=\theta_1$ 
the commutativity in \eqref{comm_h} holds. By Lemma~\ref{Lemma_imp}, this forces $\theta_1$ to be discontinuous since we have assumed that $g$ does not have the USP. When we have the product topology on $\overleftarrow{X}_U$,  $H$ is continuous if and only if all of its factor or coordinates maps are continuous. Here $\theta_1$ is discontinuous, and hence $H$ is discontinuous. $\blacksquare$.

\vspace{-0.5cm}

\section{Solution-Map as an Embedding}  \label{Sec_Stability}

When a driven system $g$ induces a causal mapping, then by Theorem~\ref{Thm2}, $g$ has the USP, and hence the solution-map $\Psi$ is well-defined. Here, we show that $\Psi$ is a continuous mapping and an embedding when $H$ is a causal mapping and a causal embedding, respectively. The continuity of $\Psi$, of course, means close-by inputs are mapped to close-by solutions -- robustness to noise of the driven system when it has the USP.  
The continuity of the solution-map $\Psi$ when $X$ is a compact subset of $\mathbb{R}^N$ is known as the fading memory property in some literature \cite{boyd1985fading,grigoryeva2018echo}.

\begin{Theorem} \label{Thm3}
Let $g$ be a compactly driven system that has a causal mapping  (see \eqref{eqn_H}). Then 
	the solution-map $\Psi: \overline{U} \to \overline{X}$ is continuous and is given by 
	\begin{equation} \label{eqn_Psi}
		\Psi(\bar{u}) = (\ldots,h(\cev{u}^{-1}),h(\cev{u}), h(\cev{u}^1),\dots),
	\end{equation}
where $h(\cev{u}^r)$ is the $r^{\text{th}}$ component of $\Psi(\bar{u})$ and $h(\cev{u})$ is the $0^{\text{th}}$ component. If further, a causal embedding $H$ exists then the solution-map $\Psi: \overline{U} \to \overline{X}$ is an embedding. 
\end{Theorem}

{\bf Proof. } Fix $\bar{u}$. By Theorem~\ref{Thm2}, $g$ has the USP since it has a causal mapping. Next, $g$ has the USP if and only if the encoding map $\mathcal{E}$ takes values in $\mathsf{S}_X$ and hence $h = \mathit{i}^{-1} \circ \mathcal{E}$ is well-defined on $\overleftarrow{U}$. Therefore, $h(\cev{u}^k)=\mathit{i}^{-1} \circ \mathcal{E}(\cev{u}^k)$ is the value of the solution at time $k$. Hence $\Psi(\bar{u})$ is as in \eqref{eqn_Psi}.  It is continuous since the coordinate functions are all continuous in view of  Lemma~\ref{Lemma_imp}.

To show that $\Psi$ is an embedding, we prove that $\Psi: \overline{U} \to \overline{X}_U$ is a homeomorphism. First, we shall show that $\Psi$ is injective. Suppose 
$\Psi(\bar{u}) = \Psi(\bar{v})$ for some $\bar{u} \not= \bar{v}$. Since $\bar{u} \not= \bar{v}$ there exists an $k \in \mathbb{Z}$ such that $u_k\not=v_k$. This would imply $H(\cev{u}^k) \not= H(\cev{v}^k)$ and this is a contradiction since $H: \overleftarrow{U} \to \overleftarrow{X}_U$ is a homeomorphism as it is a causal embedding.  Hence $\Psi$ is injective. The causal embedding $H : \overleftarrow{U} \to \overleftarrow{X}_U$ is continuous and surjective. Note that the solution-map $\Psi : \overline{U} \to \overline{X}$ is related to $H$ by $\pi_k(\Psi(\bar{u})) = \pi_1(H(\cev{u}^k))$.  Since $H$ and $\pi_1$ are continuous and surjective, the coordinate maps of $\Psi$ are all both continuous and surjective. 
Hence $\Psi: \overline{U} \to \overline{X}_U$ is continuous and surjective as well.  Since a continuous bijection between Hausdorff spaces is a homeomorphism, $\Psi:\overline{U} \to \overline{X}_U$ is a homeomorphism. $\blacksquare$

We next present a result to bring out another significant role of the USP by examining a stability notion that can be defined even when $g$ does not have the USP. Here, we are interested not in how a single-valued mapping like $h$ or $H$ behaves but rather in how the collection of solutions behave as a function of the input. 
More precisely, where multiple solutions can emerge for an input $\bar{u}$, we define the set-valued map $\bar{u} \mapsto 
 \{X_n(\bar{u})\}$. By Lemma~\ref{Lemma_soln}, $x\in X_k(\bar{u})$ if and only if there exists a solution $\{x_n\}$ for the input $\bar{u}$ such that $x_k =x$.  We would be interested in examining the continuity of the mapping 
 $\bar{u} \mapsto 
 \{X_n(\bar{u})\}$ when $g$ does not have the USP. The $n^{\text{th}}$ coordinate of this mapping is given by  $\mathcal{E}(\cev{u}^n)$. Hence it is sufficient to consider the continuity of the mapping $\cev{u}^k \mapsto \mathcal{E}(\cev{u}^k)$ when the product topology is employed. So the continuity of $\mathcal{E}$ reflects on whether small changes in $\cev{u}$ would result in a `proportionate' response in the ensemble of solution values and not a single solution -- a notion of input related stability as defined next in Definition~\ref{Def_Girst}.

\begin{Definition} \rm \label{Def_Girst} Let $g$ be a driven system and $\mathcal{E}$ be the encoding map. If $\mathcal{E}: \overleftarrow{U} \to \mathsf{H}_X$ is continuous ($\mathsf{H}_X$ is equipped with the topology generated by the Hausdorff metric) then we say $g$ has the global-input related stability (GIRST).
\end{Definition}

The local continuity of $\mathcal{E}$ is addressed in \cite{Manju_2020}. Here, we consider the global continuity of $\mathcal{E}$ and show that the USP of $g$ is equivalent to $g$ having the GIRST for a class of systems.  The proof of Proposition~\ref{Prop_GIRST} is an adoption of the result \cite[Theorem 1]{Manju_2020} but we are able to show that with fewer hypotheses on $g$. 
We also note that we do not need compactness on $U$ to prove the result.

\begin{Proposition} \rm \label{Prop_GIRST}
Consider a driven system $g$. Suppose $g$ has the USP, then $g$ has the 
GIRST. Conversely, if $g$ has the GIRST and is such that it also has precisely one solution for at least one input $\bar{u}\in \overline{U}$, then $G$ has the USP. In particular, if $g$ has precisely one solution for some input $\bar{u}$, then the GIRST and the USP are equivalent. 
\end{Proposition}  

\vspace{-0.5cm}
{\bf Proof. }  (USP $\Longrightarrow$ GIRST). 
We know that the map $\mathcal{E}: \overleftarrow{U} \to \mathsf{H}_X$ is continuous when $g$ has the USP since $\mathcal{E}$ is upper-semicontinuous (by Lemma~\ref{Lemma_USC2}) and taking values in $\mathsf{S}_X$ (by Remark~\ref{Rem_cont}). 

(GIRST $\Longrightarrow$ USP). When $g$ has the GIRST and has exactly one solution for some input $\bar{v}$, then fix $\bar{v}$. Assume that there exists an input $\bar{u}$ that has two different solutions. That is there exists an $k\in \mathbb{Z}$ such that
 $\mathcal{E}(\cev{u}^k)$ contains at least two elements of $X$. Without loss of generality, let $k=0$. Hence the encoding map $\mathcal{E}$ evaluated at $\cev{u}$ is bounded away from all singleton subsets of $X$, i.e.,
 $r:=\inf_{y \in \mathsf{S}_X} d_H(\mathcal{E}(\cev{u}), y)>0$. 
 
We next define a family of left-infinite sequences
$\{\cev{w}_n\}_{n\in \mathbb{N}}$ using $\bar{v}$, where
 each left-infinite sequence in the family
$\cev{w}_n : = (\ldots, v_{-2},v_{-1},u_{-{n}},\ldots,u_{-2},u_{-1})$. Since, the first $n$ elements of $\cev{w}_n$ are identical to that of $\cev{u}$, $\cev{w}_n \to v$ since $\overleftarrow{U}$ is equipped with the product topology. Let $\bar{w}_n$ and $\bar{v}$ be any fixed choice of the right-infinite extension of $\cev{w}_n$ and $\cev{v}$ respectively so that the time-indices are preserved, i.e.,  the $-1^{\text{th}}$ element of $\bar{v}$ is $v_{-1}$ and the $(-n-1)^{\text{th}}$ element of $\bar{w}_n$ is $v_{-1}$. Thus,
\begin{equation} \label{Eqn_Girst}
	\mathcal{E}(\cev{w}_n^{-n}) = \bigcap_{j>0} \phi_{\bar{w}_n}(-n,-j,X) = \bigcap_{j>0} \phi_{\bar{v}}(0,-j,X) = \mathcal{E}(\cev{v}),
\end{equation}
By our hypothesis, $\mathcal{E}(\cev{v}) \in \mathsf{S}_X$, so in view of \eqref{Eqn_Girst},  $\mathcal{E}(\cev{w}_n^{-n}) \in \mathsf{S}_X$ for all $n\in \mathbb{N}$. When $\mathcal{E}(\cev{w}_n)\in \mathsf{S}_X$  it follows from the definition of the positive real $r$ above that $d_H(\mathcal{E}(\cev{u}),\mathcal{E}(\cev{w}_n)) \ge r> 0$ for all $n$. This implies $\mathcal{E}_n(\cev{w}_n)$ does not converge to $\mathcal{E}(\cev{u})$ although we had  $\cev{w}_n \to \cev{u}$. This contradicts our assumption  that $\mathcal{E}(\cev{u})$ is continuous at $\cev{u}$, and hence also contradicts the hypothesis that $g$ has the GIRST. 
$\blacksquare$

\vspace{-1cm}

\section{Conclusions} \label{Sec_Conclusions}

Often while collecting data from natural systems and in other examples in the human-made world,  one does not have access to the full states of a complicated system. If one has to forecast data from such systems, one could model the unobserved states as stochastic quantities. Alternatively, if the data is mapped onto a higher-dimensional space through a time-delay observation map, the unobserved coordinates of the data source may get revealed. For instance, if the observed data is a partial scalar observable of motion along an attractor of a smooth dynamical system, Takens embedding could hold. Discrete-time state space models that map input data onto a higher dimensional space can also embed attractors \cite{grigoryeva2021learning} and also help in information processing tasks that involve more general inputs, and have the state-of-the-art status for many applications. At the heart of every discrete-time state space model is a driven dynamical system.  In this paper, we have shown that if the input originates from a compact space, and if every input entails exactly one solution, i.e., if the driven system has the unique solution property,  
then it is equivalent to representing information without distortion, which we have formalised through the existence of the causal mapping.
 In particular, if the driven system has the unique solution property and is also invertible for each fixed state-variable, it induces an infinite-delay observable through which one can obtain the coordinate mappings of a function (a causal embedding) that embeds sequential data contained in a compact input space in the solution space of driven system's state space.  Like the essence of a conjugacy in Takens embedding, the causal embedding ensures that the dynamics on the nonautonomous attractor in the driven system's state space is topologically conjugate to that in the corresponding input sequence.  The possibility of a causal embedding sheds a new and sharp light on the question of whether one could move towards temporal data processing without loss of information in the field of reservoir computing.
 When this happens, we believe the field would be put on a more appropriate foundation than it had before.
 
Also, the results are general, and any driven dynamical system with specific properties can give rise to a causal mapping. This also answers why random choices of driven systems actually work in information processing methods in \cite{jaeger2004harnessing,maass2002real}, and several thousands of other publications that have followed.

 On a separate front, when a driven system models a continuous perturbation of an irreducible attracting set of an autonomous dynamical system, we provide sufficient and necessary conditions for a uniform attractor to emerge. In a future work, we intend to extend all the results in this paper to the continuous time-setting when actually the product spaces are not metrizable.

 \vspace{-0.5cm}
{\bf Acknowledgements} The author sincerely thanks the referees for several suggestions on improving the presentation in the paper. The research was supported by the National Research Foundation of South Africa under Grant UID 115940.

 \vspace{-0.5cm}

\bibliographystyle{pnas-new}
\bibliography{pnas}

\end{document}